\documentclass[a4paper,12pt]{article}
\usepackage{latexsym,amsfonts,amsthm,amsmath,amscd,amssymb,multirow}
\newtheorem{lemma}{Lemma}[section]

\newtheorem{prop}[lemma]{Proposition}

\newtheorem{defn}[lemma]{Definition}

\def\C{{\mathbb C}}
\def\R{{\mathbb R}}

\def\N{{\mathbb N}}

\def\qgm{{q^{ \Gamma/2}}}
\def\qmgm{{q^{ -\Gamma/2}}}
\def\qpmgm{{q^{ \pm\Gamma/2}}}
\def\qmpgm{{q^{ \mp\Gamma/2}}}
\def\qg{{q^{ \Gamma}}}
\def\qmg{{q^{ -\Gamma}}}
\def\opo{{\cal O}}
\def\perm{{\cal S}}

\title{Braid group representations from a deformation of the harmonic oscillator algebra}
\author{M. Tarlini\footnote{\small INFN Sezione di Firenze, 
Dipartimento di Fisica e Astronomia, Universit\`a degli Studi Firenze, email: Marco.Tarlini@fi.infn.it}}

\begin{document}

\maketitle

\begin{abstract}
We describe a new technique to obtain representations of the braid group $B_n$ from the $\Re$--matrix of a quantum deformed algebra of the one 
dimensional harmonic oscillator.  We consider the action of the $\Re$--matrix not on the tensor product of representations of the algebra, that in the 
harmonic oscillator case are infinite dimensional, but on the subspace of the tensor product corresponding to the lowest weight vectors. 

\end{abstract}

\thispagestyle{empty}

\section{Introduction}

Given a vector space and a matrix acting on it a generator of the Artin braid group $B_n$ \cite{A} can be represented on the $n$-tensor product of the vector space  
if and only if the matrix satisfies the Yang-Baxter equation, the original main reason for the quantum groups was to solve this equation \cite{FRT}, and to study 
links invariants \cite{KR}, \cite{RT}.

In this paper we want to give a new way to find representation of the braid group $B_n$ using a quantum group version \cite{h1q} of the harmonic oscillator algebra. 
This Lie algebra has four generators; we consider, besides the ladder operators, the constant operator and the Hamiltonian as generators. The irreducible representations 
of this algebra are infinite dimensional and it is not possible to use the standard methods based on the quantum $\Re$--matrix formalism without a regularization procedure,
for instance a connection between the quantum deformed harmonic oscillator and some links invariants is given in \cite{GS}. 

Here we use a different approach which takes into account the fact that lowest (highest) weight spaces (they are finite dimensional) in the $n$-tensor product of representations 
of the quantum algebra are $B_n$ invariant and they define a representation of the braid group. 

This idea has been followed in \cite{JK} for the $U_q(sl_2)$ quantum algebra, where the authors indicate with ${\bf W}_{n,l}$ the highest weight spaces and prove the isomorphism
between ${\bf W}_{n,1}$ and the reduced Burau representation \cite{Burau},\cite{KT} of $B_n$ and between ${\bf W}_{n,2}$ and the Lawrence-Krammer-Bigelow representation 
\cite{LKB}, \cite{KT}. 

The advantage of working with the harmonic oscillator algebra is that we are able to construct the $B_n$ representation spaces analogous to the lowest weight spaces ${\bf W}_{n,l}$
very easily. Indeed it is possible to start with different representations of the quantum algebra in each tensor space, and to reach lowest weights space of higher level 
(corresponding to their parameter $l>2$). We recover in the simplest case, and with all the representations in the tensor product being equal each others, the reduced Burau 
representation and in a second step the Lawrence-Krammer-Bigelow representation but with a parameter less. 

In general following the formula (\ref{sigmaopo}) it is possible to quantize a general classical representations of $B_n$. 

The paper is organized as follows: in Section 2 we recall the quantum deformation of the harmonic oscillator algebra and its irreducible representation, in Section 3 
we write explicitly the $\Re$--matrix and the assignment of the generators $\sigma_i$ of the braid group $B_n$. In Section 4  we describe the lowest weight vectors 
of the quantum algebra, and in Section 5 we introduce the operator $\opo$ and the reduction of the representations on the $n$-tensor product. In Section 6 we 
write the general formula (\ref{sigmaopo}) that allows to build the representations of the braid group. Finally in Section 7 we give some explicit examples.

\section{Quantum deformed algebra and his irreducible representations}\label{qda_irreps}

In \cite{3dc,h1q} we obtained a quantum group deformation of the Lie algebra of the harmonic oscillator by a contraction of the quantum algebra $sl(2)_q$.
Here we recall the structure of Hopf algebra of this deformation that we write in the rational form and call $ho_q$.

The associative algebra $ho_q$ is defined over $\C$ with generators $\alpha^\pm\,,\epsilon\,,\qgm\,$, $\qmgm$ and $1$ with the relations 
$\qpmgm \qmpgm=1$ and with the following commutators, where $q$ is a nonzero parameter:
\begin{equation}\label{hoq}
[\alpha^-,\alpha^+]=\frac{\qg-\qmg}{q-q^{-1}}\equiv[\Gamma]_q\ ,\quad\quad [\epsilon,\alpha^{\pm}]=\pm\,\alpha^\pm\ ,\quad\quad [\qpmgm, \cdot]=0\ .
\end{equation}

The Hopf algebra coproduct is given by 
\begin{gather}\label{hopf}
 \Delta\, \alpha^\pm=\alpha^\pm\otimes\, \qgm + \qmgm\, \otimes\alpha^\pm\ ,\\
 \Delta\, \epsilon =\epsilon\,\otimes 1 +1\otimes\, \epsilon\ , \quad
 \Delta\, \qpmgm =\qpmgm\otimes \qpmgm\ ,
 \end{gather}
 with counit of $\alpha^\pm$ and $\epsilon$ equal to zero and counit of $\qpmgm$ equal to one, the antipod is
 \begin{gather}
 S(\alpha^\pm)=-\alpha^\pm\ , \quad S(\qpmgm)=\qmpgm , \quad S(\epsilon)=-\epsilon\ .
\end{gather}

The algebra $ho_q$ has two Casimirs, one is $\qgm$ the other is
\begin{equation}
 C_q=[\Gamma]_q\; \epsilon -\alpha^+\,\alpha^-\ .
\end{equation}

The representations of the braid groups will be derived starting from the representations of this algebra.
The irreps are labeled by two numbers: $\gamma$ and $c$, we call ${\bf H}^{(\gamma,c)}$ 
the module freely generated by a set of vectors $h^{(\gamma,c)}_m\,,\ m\in \N_0$.

The representation of $ho_q$ on ${\bf H}^{(\gamma,c)}$ is then
\begin{eqnarray}\label{irreps}
 \alpha^-\cdot h^{(\gamma,c)}_m &=& [\gamma]_q^{1/2}\,m^{1/2}\, h^{(\gamma,c)}_{m-1}\ ,\nonumber\\ 
 \alpha^+\cdot h^{(\gamma,c)}_m &=& [\gamma]_q^{1/2}\,(m+1)^{1/2}\, h^{(\gamma,c)}_{m+1}\ ,\\
 \epsilon\cdot h^{(\gamma,c)}_m &=&(m+c)\ h^{(\gamma,c)}_m\ . \nonumber
\end{eqnarray}
With $\qpmgm\cdot h^{(\gamma,c)}_m= q^{\pm\gamma/2}\ h^{(\gamma,c)}_m$  and  $C_q\cdot h^{(\gamma,c)}_m=[\gamma]_q\,c\, h^{(\gamma,c)}_m$. 

The generators of this representation can be obtained from the vector $h^{(\gamma,c)}_0$ by the action of $\alpha^+$, namely 
$h^{(\gamma,c)}_m=([\gamma]_q^m m!)^{-1/2}(\alpha^+)^m\, h^{(\gamma,c)}_0$.

We assume $\gamma,c\in \R$, we choose $q\in \R$ ($0<q<1$), then this representation is hermitian with respect to the star involution 
$(\alpha^\pm)^*=\alpha^\mp\ ,(\epsilon)^*=\epsilon$ and $(\gamma)^*=\gamma$. 
The star involution is an anti--homomorphism and it fulfills $\Delta \circ *=(*\otimes\,*)\circ \Delta$ and $S\circ *=*\circ S^{-1}$. The scalar product is given by
$<h^{(\gamma,c)}_m|h^{(\gamma,c)}_{m'}>=\delta_{m\,m'}$.

\section{The $\Re$--matrix}\label{Rm}

In \cite{h1q} the quasitriangular $\Re$--matrix is obtained as a contraction limit of the $\Re$--matrix of $sl_q(2)$,  the result is
\begin{equation}\label{rmatrix}
\Re=q^{ -({\textstyle \epsilon}\,\otimes\,\Gamma + \Gamma\,\otimes\,{\textstyle \epsilon})} \exp[(q-q^{-1})\,(\qgm\,\otimes\,\qmgm)\,\alpha^-\,\otimes\,\alpha^+]\ .
\end{equation}
It is a general feature that the $\Re$--matrix it is not an element of the tensor product of the quantum algebra given in the rational form, nevertheless it acts on
any tensor product of algebra representations (see for instance Chapter 9 and 10 of the book \cite{CP} for a general reference).
We are interested in the action of $\Re$ on the representations given in (\ref{irreps}), this leads us to consider 
$\Gamma\cdot h^{(\gamma,c)}_m= \gamma\ h^{(\gamma,c)}_m$. 

$\Re$ is then an endomorphism on ${\bf H}^{(\gamma_1,c_1)}\,\otimes\, {\bf H}^{(\gamma_2,c_2)}$. If ${\cal P}$ is the permutation on the tensor product, 
the action of ${\cal P}\, \Re$ gives an action of a generator of the braid group. 

We write here the action of ${\cal P}\Re$ on an element $h^{(1)}_m\,\otimes\,h^{(2)}_{m'}\in{\bf H}^{(\gamma_1,c_1)}\,\otimes\, {\bf H}^{(\gamma_2,c_2)}$ 
, to simplify the notations we use $h_m^{(i)}$ for $h_m^{(\gamma_i,c_i)}$ :
\begin{equation}\label{p_er}
\begin{split}
 {\cal P}\Re\;\cdot h^{(1)}_m\,\otimes\,h^{(2)}_{m'}&= q^{-((m+c_1)\gamma_2 +(m'+c_2)\gamma_1)} \sum^m_{k=0}{\binom {m+k-1}{m-1}}^{1/2} {\binom {m'+k}{m'}}^{1/2}\cdot \\
&\cdot  (1-q^{-2\gamma_1})^{k/2} (q^{2\gamma_2}-1)^{k/2}\,h^{(2)}_{m'+k}\,\otimes\,h^{(1)}_{m-k}\ . 
\end{split}
\end{equation}

The Artin braid group $B_n$ is defined as the group generated by the $n-1$ generators $\sigma_1,\sigma_2,\cdots\sigma_{n-1}$ and the {\it braid relations}: $\sigma_i\sigma_j=
\sigma_j\sigma_i$ for $|i-j|>1$ and $\sigma_i\sigma_{i+i}\sigma_i=\sigma_{i+1}\sigma_i\sigma_{i+1}$ for $i=1,\cdots, n-2$ (see Section 1.1 of \cite{KT}).

The Yang--Baxter equation $\Re_{12}\Re_{13}\Re_{23}=\Re_{23}\Re_{13}\Re_{12}$ provides the braid relations of the generators $\sigma_i$ of the $B_n$ braid group 
presentation, then a representation of $B_n$ on ${\bf H}^{\otimes\,n}$ is given by the following identification
\begin{equation}\label{sigma}
 \sigma_j=1^{\otimes\,j-1}\,\otimes\,{\cal P}\,\Re\,\otimes\,1^{\otimes\,n-j-1}\ .
\end{equation}

The inverse of the $\Re$--matrix and consequently the inverse of $\sigma_i$ are obtained by the expression $\Re^{-1}=(S\,\otimes 1)\,\Re$; from (\ref{rmatrix}) 
it results in the exchange $q\rightarrow q^{-1}$ in $\Re$ and then in $\sigma_i$.
In the following section we describe a method to find a finite dimensional space of representation of $B_n$ starting from (\ref{sigma}) and the infinite 
dimensional representations (\ref{irreps}) of $ho_q$. The main feature is the fact that the operator ${\cal P}\,\Re$, and hence also the maps $\sigma_i$, commute with 
the action of the $ho_q$ algebra on the tensor product, this action is given by the coproduct.

\section{Lowest weight spaces of $ho_q$}\label{sec_weight_spaces}

The action of the generators of $ho_q$ on the $n$--tensor product of representations is given by the iterated coproduct defined as $\Delta^{(n)} = 
(\Delta^{(n-1)}\,\otimes  1) \Delta$ with $n\geq 2$ and $\Delta^{(2)}\equiv\Delta$.

\begin{defn}\label{weight_space}
The weight space ${\bf W}^{(\gamma,e)}$ corresponding to the weights $\gamma$ and $e$ is defined as ${\rm ker}(\Delta^{(n)}\Gamma -\gamma)\cap 
{\rm ker}(\Delta^{(n)}\epsilon -e)\subset 
\bigoplus_{\perm\in S_n}\bigotimes_{i=1}^n {\bf H}^{(\gamma_{\perm(i)},c_{\perm(i)})}$, where $\perm\in S_n$ is a permutation of $(1,\cdots,n)$.
\end{defn}

The generator $\epsilon$ and $\Gamma$ are primitive then ${\bf W}^{(\gamma,e)}$ is the span of vectors of the form 
\begin{equation}\label{vector_n}
h^{\perm(1)}_{m_1}\,\otimes h^{\perm(2)}_{m_2}\cdots\otimes\,h^{\perm(n)}_{m_n}\quad  \text{with} \quad \left\{
\aligned  \gamma_1+\cdots+\gamma_n&=\gamma\ , \\m_1+\cdots +m_n&=e-c_1-\cdots -c_n\ .
\endaligned \right.
\end{equation}
Here we used the short notation $h^{(i)}_m$ for $h^{(\gamma_i,c_i)}_m$.

Then we have that $e\geq \sum_i^n c_i$ and $e - \sum_i^n c_i$ is a nonnegative integer, moreover the value $\gamma$ is fixed by the $\gamma_i$. 

In the following we fix the representations ${\bf H}^{(\gamma_i,c_i)}$.

Note that the Casimir $C_q$ is not primitive therefore a generic vector in ${\bf W}^{(\gamma,e)}$ is not an eigenvector of $\Delta^{(n)} C_q$, in particular of the term 
$\Delta^{(n)}\alpha^+\,\Delta^{(n)}\alpha^-$.  

\begin{defn}\label{lowest_weight}
The lowest weight space corresponding to $e$ is the space ${\bf V}^{(e)}\subset {\bf W}^{(\gamma,e)}$ of vectors $v_0^{(e)}$ such that 
$\Delta^{(n)} \alpha^- \cdot v_0^{(e)}=0$. 
\end{defn}

For every vector $v_0^{(e)}\in{\bf V}^{(e)}$ we have that 
\begin{equation}\label{deltaC}
 \Delta^{(n)} C_q\cdot\, v_0^{(e)}= (\frac{\Delta^{(n)}\, q^\Gamma-\Delta^{(n)}\, q^{-\Gamma}}{q-q^{-1}}\; 
 \Delta^{(n)}\epsilon -\Delta^{(n)}\alpha^+\,\Delta^{(n)}\alpha^-)\cdot\,v_0^{(e)}=[\gamma]_q\,e\;v_0^{(e)}\ .
\end{equation}
From now on we will indicate the eigenvalue of $\Delta^{(n)} C_q$ neglecting the $[\gamma]_q$ factor, namely we consider $e$ as the eigenvalue of $v_0^{(e)}$.
Starting from a vector $v_0^{(e)}\in {\bf V}^{(e)}$ one can built a representation of $ho_q$ by applying the operator $\Delta^{(n)}\alpha^+$, namely we define 
\begin{equation}\label{vi}
v_m^{(e)}=([\gamma]^m_q\, m!)^{-1/2}(\Delta^{(n}\alpha^+)^m \, v_0^{(e)}\ ,
\end{equation}
we notice that $v_m^{(e)}\in {\bf W}^{(\gamma,e+m)}$ is not an element of ${\bf V}^{(e+m)}$, indeed 
$\Delta^{(n)} \alpha^-\cdot\; v_m^{(e)}=[\gamma]_q^{1/2}\,m^{1/2}\, v_{m-1}^{(e)}\neq0$.
The eigenvalue of the Casimir on $v_m^{(e)}$ is $e$ .

The relations (\ref{irreps}) are fulfilled with the substitutions: $c\rightarrow e$, $h^{(c)}_m\rightarrow v_m^{(e)}$, 
$\alpha^+\rightarrow \Delta^{(n)} \alpha^+$, $\alpha^-\rightarrow \Delta^{(n)} \alpha^-$, $\epsilon\rightarrow \Delta^{(n)} \epsilon$ and the scalar product is
$<v_m^{(e)}|v_{m'}^{(e)}>=\delta_{m\,m'}$. 

\bigskip
{\bf Examples:} For $n=2$ the lowest value for $e$ is $c_1+c_2$. The lowest weight space ${\bf V}^{(c_1+c_2)}$ is generated by $v_0^{(c_1+c_2)}=h_0^{(1)}\otimes h_0^{(2)}$ 
(and the permutation $h_0^{(2)}\otimes h_0^{(1)}$), namely $\Delta\epsilon\cdot\, v_0^{(c_1+c_2)}=(c_1+c_2)\; v_0^{(c_1+c_2)}$ and  
$\Delta\alpha^-\cdot\, v_0^{(c_1+c_2)}=0$.

The next value for $e$ is $c_1+c_2+1$ then we write $v_0^{(c_1+c_2+1)}=a\;h_0^{(1)}\otimes h_1^{(2)} + b\;h_1^{(1)}\otimes h_0^{(2)}$ so that 
$\Delta\epsilon\cdot\, v_0^{(c_1+c_2+1)}=(c_1+c_2+1)\; v_0^{(c_1+c_2+1)}$. Imposing the lowest weight relation we get $v_0^{(c_1+c_2+1)}$ proportional to
$q^{\gamma_2/2}\,[\gamma_1]_q^{1/2}\;h_0^{(1)}\otimes h_1^{(2)} - q^{-\gamma_1/2}\,[\gamma_2]_q^{1/2}\;h_1^{(1)}\otimes h_2^{(2)}$. We have that ${\bf V}^{(c_1+c_2+1)}$
is generated by this $v_0^{(c_1+c_2+1)}$ and the permutation obtained by the exchange $\gamma_1 \leftrightarrow \gamma_2$ and $(1)\leftrightarrow (2)$.

\section{The operator $\opo$ and the reduction of the tensor product of representations}

We start with $n=2$, we take the representations ${\bf H}^{(\gamma_1,c_1)}$ and ${\bf H}^{(\gamma_2,c_2)}$. 
The idea is to find an operator $\opo\in ho_q\,\otimes\,ho_q$ that commutes with
$\Delta \alpha^-$ and that works as a ladder operator on $\Delta\epsilon$ so that $\opo$ maps vectors of ${\bf V}^{(e)}$ in vectors of ${\bf V}^{(e+1)}$.  
The following proposition derives from a straightforward computation:

\begin{prop}
The operator
\begin{equation}\label{opo}
 \opo=q^{-\Gamma/2}\,[\Gamma]_q^{-1/2}\alpha^+\,\otimes\,[\Gamma]_q^{1/2} - [\Gamma]_q^{1/2} \,\otimes\,q^{\Gamma/2}\,[\Gamma]_q^{-1/2}\alpha^+
\end{equation}
satisfies the following commutators:
\begin{equation}\label{commopo}
[\Delta \alpha^-,\opo]=0\;, \qquad [\Delta \epsilon,\opo]=\opo\ ,
\end{equation} 
and obviously
\begin{equation}
[\Delta \Gamma,\opo]=0\;,\qquad [\Delta \alpha^+,\opo]=0\ .
\end{equation}
\end{prop}

Following Definition \ref{lowest_weight} the space ${\bf V}^{(c_1+c_2)}$ is the lowest weight space with $e=c_1+c_2$ and $\gamma=\gamma_1+\gamma_2$, 
where $\gamma_i\,,c_i$ define the Casimirs eigenvalues on the representations ${\bf H}^{(\gamma_i,c_i)}$.

From the vector $v_0^{(c_1+c_2)}=h^{(1)}_{0}\otimes h^{(2)}_{0}\in{\bf V}^{(c_1+c_2)}$ we can obtain  $v_0^{(c_1+c_2+1)}=[\gamma]_q^{-1/2}\opo\,v_0^{(c_1+c_2)}$, 
with the eigenvalue of the Casimir $\Delta C_q$  given by $c_1+c_2+1$ and  $v_1^{(c_1+c_2)}=[\gamma]_q^{-1/2}\,\Delta\alpha^+\,v_0^{(c_1+c_2)}$, 
with the eigenvalue of the Casimir given by $c_1+c_2$. 

One can start from the vector ${v'}_0^{(c_1+c_2)}=h^{(2)}_{0}\otimes h^{(1)}_{0}\in{\bf V}^{(c_1+c_2)}$
and repeat the same construction.

The vectors $v_0^{(c_1+c_2+1)}$ and $v_1^{(c_1+c_2)}$ together with ${v'}_0^{(c_1+c_2+1)}$ and ${v'}_1^{(c_1+c_2)}$  have the eigenvalue of $\Delta\epsilon$ equal to 
$c_1+c_2+1$ and they span ${\bf W}^{(\gamma,c_1+c_2+1)}=\langle \; h_0^{(\perm(1))}\otimes h_1^{(\perm(2))},\,h_1^{(\perm(1))}\otimes h_0^{(\perm(2))}\,;\ \perm\in S_2\;\rangle$, 
namely the weight space with $e=c_1+c_2+1$.

In general  ${\bf W}^{(\gamma,c_1+c_2+N)}$ is spanned by the $2(N+1)$ vectors $h_m^{\perm(1)}\otimes h_{m'}^{\perm(2)}$ for $m+m'=N$. It can be reduced in the combination
of vectors of the representations with Casimir from $c$ to $c+N$ of the form $v_i^{(c_1+c_2+j)}$  with $i+j=N\ (i,j\geq0)$ obtained from 
$v_0^{(c_1+c_2)}=h_0^{\perm(1)}\otimes h_0^{\perm(2)}$ using (\ref{vi})  and
\begin{equation}\label{v0}
v_0^{(c_1+c_2+k)}=([\gamma]^k_q\, k!)^{-1/2}\opo^{\,k} \, v_0^{(c_1+c_2)}\ .
\end{equation}

From the commutation relation $[\opo^*,\opo]=\Delta [\Gamma]_q$ and $[\Delta\alpha^-,\opo]=0$ we can prove that these vectors verify  
\begin{equation}
 <v_i^{(c_1+c_2+j)}|v_{i'}^{(c_1+c_2+j')}>=\delta_{i\,i'}\delta_{j\,j'}\ .
\end{equation}
The rising of $c_1+c_2$ by $j$ steps using $\opo$ and the rising from $0$ to $i$ using $\Delta\alpha^+$ can be made in any order because of the commutativity of 
$\opo$ with $\Delta\alpha^+$. 

\bigskip
We have then the following proposition:
\begin{prop}\label{reduction}
For $n=2$ we define ${\bf U}^{(j,N)}\subset {\bf W}^{(\gamma,c_1+c_2 +N)}$ as the space spanned by the vectors 
$\langle\; (\Delta\alpha^+)^{(N-j)}\;\opo^{\,j}\,v_0^{(c_1+c_2)}\; \rangle$, with $v_0^{(c_1+c_2)}\in {\bf V}^{(c_1+c_2)}$. The vectors in ${\bf U}^{(j,N)}$ 
are eigenvectors of $\Delta C_q$ and $\Delta \epsilon $ with eigenvalue respectively $c_1+c_2+j$ and $c+N$.
We have that 
\begin{equation}\label{eq_reduction}
{\bf W}^{(\gamma,c_1+c_2+N)}=\bigoplus_{j=0}^N{\bf U}^{(j,N)}\ .
\end{equation}
\end{prop}

We treat now the case of a generic $n$.

One can repeat what was done for $n=2$ using the operators 
\begin{equation}
\Delta^{(n)}\alpha^+ \quad \text{and}\quad \opo_k=1^{\otimes k-1}\otimes\underbrace{\opo}_{k,k+1}\otimes\; 1^{\otimes n-k-1}\,, 
\end{equation}

$\Delta^{(n)}\alpha^-$, $\Delta^{(n)}\epsilon$ and $\opo_k$ still satisfy the commutation relations (\ref{commopo}) for every $k$.  

Fixing the Casimirs $c_i$ of the representations in the tensor product we call ${\bf c}$ the sum of them: ${\bf c}=c_1+c_2+\cdots +c_n$. 

The number ${\cal N}_{n,N}$ of vectors of the form (\ref{vector_n}) (apart the action of ${\cal S}\in S_n$) that span ${\bf W}^{(\gamma,{\bf c}+N)}$ is the number of 
different ways to add (taking into account the order) $n$ nonnegative integers $(m_1,m_2,\cdots,m_n)$ to get $N$, namely
\begin{equation}
{\cal N}_{n,N}=\binom{n+N-1}{n-1} \,,
\end{equation}
the vectors obtained by the action of degree $j$ monomials $\opo_1^{j_1}\opo_2^{j_2}\cdots \opo_{n-1}^{j_{n-1}}$ on $v_0^{({\bf c})}$ with nonnegative $j_i$ and 
$j=j_1+\cdots +j_{n-1}$ are eigenspaces of constant Casimir equal to ${\bf c}+j$. They are in a number equal to
\begin{equation}
{\cal M}_{n,j}=\binom{n+j-2}{n-2} \,.
\end{equation}
We have to act with $(\Delta^{(n)}\alpha^+)^{N-j}$ on the vectors with constant Casimir ${\bf c}+j$ to get vectors in ${\bf W}^{(\gamma,{\bf c} +N)}$, 
namely with $\Delta^{(n)}\epsilon$ equal to ${\bf c}+N$ . If we sum up all the value of $j$ from $0$ to $N$ we recover ${\cal N}_{n,N}\,$:
\begin{equation}
{\cal N}_{n,N}=\sum_{j=0}^N {\cal M}_{n,j}\ .
\end{equation}

For example let us built the vectors of ${\bf W}^{(\gamma,{\bf c}+3)}$ for $n=3$ and ${\bf c}=c_1+c_2+c_3$, (they have then the eigenvalue of $\Delta^{(3)}\epsilon$ equal to 
${\bf c}+3$): 
\bigskip

\begin{tabular}{rcc}

 States                                             & Casimir  &            ${\cal N}_{3,3}=10$\\
                                                    &          &      \\
$(\Delta^{(3)}\alpha^+)^3\,v_0^{({\bf c})}$               &  ${\bf c}+0$   &                      ${\cal M}_{3,0}=1$\\       
$(\Delta^{(3)}\alpha^+)^2\,\opo_1\,v_0^{({\bf c})}$       &  ${\bf c}+1$   &      \multirow{2}{*}{${\cal M}_{3,1}=2$}\\
$(\Delta^{(3)}\alpha^+)^2\,\opo_2\,v_0^{({\bf c})}$       &  ${\bf c}+1$   &      \\
$\Delta^{(3)}\alpha^+\,\opo_1\,\opo_2\,v_0^{({\bf c})}$   &  ${\bf c}+2$   &      \multirow{3}{*}{${\cal M}_{3,2}=3$}\\
$\Delta^{(3)}\alpha^+\,(\opo_1)^2\,v_0^{({\bf c})}$       &  ${\bf c}+2$   &      \\
$\Delta^{(3)}\alpha^+\,(\opo_2)^2\,v_0^{({\bf c})}$       &  ${\bf c}+2$   &      \\
$(\opo_1)^3\,v_0^{({\bf c})}$                         &  ${\bf c}+3$   &      \multirow{4}{*}{${\cal M}_{3,3}=4$}\\
$(\opo_1)^2\,\opo_2\,v_0^{({\bf c})}$                 &  ${\bf c}+3$   &      \\
$\opo_1\,(\opo_2)^2\,v_0^{({\bf c})}$                 &  ${\bf c}+3$   &      \\
$(\opo_2)^3\,v_0^{({\bf c})}$                         &  ${\bf c}+3$   &      \hspace{3cm} .\\
\end{tabular}

The following proposition is a generalization of Proposition \ref{reduction}:
\begin{prop}\label{reduction_n}
 For generic $n$ we define ${\bf U}^{(j,N)}\subset {\bf W}^{(j,{\bf c}+N)}$ with $j\leq N$ as:
 \begin{equation}
  {\bf U}^{(j,N)}=\langle \; (\Delta^{(n)} \alpha^+)^{N-j}\,\opo_1^{j_1} \opo_2^{j_2}\cdots\opo_{n-1}^{j_{n-1}}\,v_0^{({\bf c})}\, ;\, j_1+\cdots+j_{n-1}=N\; \rangle\ ,
 \end{equation}
with $j_i\in \N_0$, where ${\bf c}=c_1+c_2+\cdots+c_n$ and $v_0^{({\bf c})}\in {\bf V}^{({\bf c})}$. 

The vectors in ${\bf U}^{(j,N)}$ are eigenvectors of $\Delta^{(n)} C_q$ and $\Delta^{(n)} \epsilon$ with eigenvalues respectively ${\bf c}+j$
and ${\bf c}+N$.
We have that:
\begin{equation}\label{eq_reduction_n}
{\bf W}^{(\gamma,{\bf c}+N)}=\bigoplus_{j=0}^N{\bf U}^{(j,N)}\ .
\end{equation}
\end{prop}

{\it Proof}:
We have that:
\begin{equation}\label{recursion}
{\bf W}^{(j,{\bf c}+N)}=\langle \; (1^{\otimes i-1}\otimes \alpha^+\otimes1^{\otimes n-i})\ {\bf W}^{(\gamma,{\bf c}+N-1)}\,;i=1,\cdots,n\;\rangle\ .
\end{equation}
Each of the $n$ vectors $(1^{\otimes i-1}\otimes \alpha^+\otimes1^{\otimes n-i})\, w^{({\bf c}+N-1)}$, where $w^{({\bf c}+N-1)}\in {\bf W}^{(\gamma,{\bf c}+N-1)}$ 
has the form (\ref{vector_n}) with $m_1+m_2+\cdots+m_n=N-1$, can be written as a linear combination of the vector $\Delta^{(n)} \alpha^+\, w^{({\bf c}+N-1)}$ and 
the $n-1$ vectors $\opo_k\,w^{({\bf c}+N-1)}$, note that the terms containing $\Gamma$ take a numerical value on vectors of the form (\ref{vector_n}). 
This implies that if the equation (\ref{eq_reduction_n}) is valid for $N-1$ we have that: 
\begin{equation}\label{induction_U}
 {\bf U}^{(j,N)}=\langle \;\Delta^{(n)}\alpha^+\, {\bf U}^{(j,N-1)}\ ,\ \opo_k\,{\bf U}^{(j-1,N-1)}\,;\, k=1,\cdots,n-1\;\rangle.
\end{equation}
With an induction procedure making use of (\ref{induction_U}), we derive that ${\bf U}^{(j,N)}$ coincides with the definition in the proposition. \qed

\bigskip
The following Lemma derives directly from the previous Proposition:
\begin{lemma}\label{lwso}
 The lowest weight spaces ${\bf V}^{({\bf c}+N)}$ are obtained as the vector spaces spanned by 
 $\langle\;\opo_1^{j_1} \opo_2^{j_2}\cdots\opo_{n-1}^{j_{n-1}}\,{\bf V}^{({\bf c})}\,;\, j_1+\cdots+j_{n-1}=N\;\rangle$,  
with $j_i\in\N_0$.
 \end{lemma}

\section{Representations of the braid group $B_n$}\label{repBn}

We are ready to build the representations of the braid group $B_n$. We use the presentation of the generators of the braid group given in (\ref{sigma}). 
From the fact that $\sigma_i$ commute with $\Delta^{(n)} \epsilon$ we derive that ${\bf W}^{(\gamma,{\bf c}+N)}$ defined in Definition \ref{weight_space}, where ${\bf c}=c_1+c_2+\cdots+c_n$ 
is the sum of the Casimir labels of the representations, is a representation space for $B_n$, moreover $\sigma_i$ commute with $\Delta^{(n)} C_q$ then this space is reducible.
The lowest weight spaces ${\bf V}^{({\bf c}+N)}$ are again $B_n$ invariant.

\bigskip
We present here the main proposition:

\begin{prop}\label{representations}
We denote ${\bf V}^{({\bf c}+N)}$ the lowest weight spaces as in the Lemma \ref{lwso}.

For each choice of the representations ${\bf H}^{(\gamma_i,c_i)}$, we obtain a representation of the braid group $B_n$ given by an automorphism of ${\bf V}^{(c+N)}$.
\end{prop}

{\it Proof:}
The building block is the computation of the  conjugation of $\opo_k$ by $\sigma_i$ given by (\ref{sigma}), for $k=i$ and $k=i\pm 1$, we recall that $\sigma_i^{-1}$ 
is obtained from $\sigma_i$ by the exchange $q\leftrightarrow q^{-1}$ : 
\begin{equation}\label{sigmaopo}
\begin{split}
\sigma_i\; \opo_i&=-(1^{\otimes\,i-1}\otimes q^{-\Gamma}\otimes q^{-\Gamma}\otimes\, 1^{\otimes\,n-i-1})\;\opo_i\;\sigma_i\quad  ,\\
\sigma_i\;\opo_{i+1}&=(1^{\otimes\,i}\otimes q^{-\Gamma}[\Gamma]_q^{-1/2}\otimes\, [\Gamma]_q^{1/2}\otimes\, 
                      1^{\otimes\,n-i-2})\;\opo_i\;\sigma_i\,+ \\
                      & +(1^{\otimes\,i-1}\otimes [\Gamma]_q^{1/2}\otimes\, [\Gamma]_q^{-1/2}\otimes\, 1^{\otimes\,n-i-1})\;\opo_{i+1}\;\sigma_i\quad  ,\\
\sigma_i\;\opo_{i-1}&=(1^{\otimes\,i-1}\otimes [\Gamma]_q^{-1/2}\otimes\, [\Gamma]_q^{1/2}\otimes\, 
                      1^{\otimes\,n-i-1})\;\opo_{i-1}\;\sigma_i\,+ \\
                      & +(1^{\otimes\,i-2}\otimes\,[\Gamma]_q^{1/2}\otimes\, q^{-\Gamma}[\Gamma]_q^{-1/2}\otimes\, 
                      1^{\otimes\,n-i})\;\opo_i\;\sigma_i\quad .
\end{split}
\end{equation}

It is clear that $\sigma_i\;\opo_j=\opo_j\,\sigma_i$ for $|i-j|>1$. 

If we act with $\sigma_i$ on the lowest weight vector $v_0^{({\bf c})}=h_0^{(1)}\otimes\, \cdots \otimes\, h_0^{(n)}\in {\bf V}^{({\bf c})}$ from (\ref{p_er}) and (\ref{sigma}) 
we get the permuted vector up a numerical factor:
\begin{equation}\label{sigma0}
\begin{split}
 \sigma_i\;&h_0^{(1)}\otimes\, \cdots\otimes\,h_0^{(i)}\otimes h_0^{(i+1)}\cdots \otimes\, h_0^{(n)}= \\
 &q^{-c_i\gamma_{i+1}-c_{i+1}\gamma_i}\, h_0^{(1)}\otimes\, \cdots\otimes\,h_0^{(i+1)}\otimes h_0^{(i)}\cdots \otimes\, h_0^{(n)}\ .
\end{split}
 \end{equation}
 
Given a vector $o^{(N)}= \opo_1^{j_1} \opo_2^{j_2}\cdots\opo_{n-1}^{j_{n-1}}\,v_0^{({\bf c})}\in {\bf V}^{({\bf c}+N)}$ we compute $\sigma_i\, o^{(N)}$ by using 
(\ref{sigmaopo}) repeatedly on $\sigma_i\, \opo_{i-1}^{j_{i-1}}\, \opo_i^{j_i}\,\opo_{i+1}^{j_{i+1}}$ to write $\sigma_i$ on the right. With $\sigma_i$ 
fully on the right we use (\ref{sigma0}) and we apply on the left side of equation (\ref{sigma0}) the combintion of operators 
$\opo_1^{k_1} \opo_2^{k_2}\cdots\opo_{n-1}^{k_{n-1}}$ that have come out from the use of (\ref{sigmaopo}).
Note that the action of $\sigma_i$ does not change the sum: $j_1+j_2+\cdots j_{n-1}=k_1+k_2+\cdots + k_{n-1}$, then we remain with vectors in ${\bf V}^{({\bf c}+N)}$. \qed

\bigskip
\section{Examples}
Here would like to present some examples of representations of $B_n$ obtained using the Proposition \ref{representations}. 
\subsection{N=1}

Let us choose as a first example $N=1$ and all $\gamma_i$'s and $c_i$'s equal. 

The representation has dimension  ${\cal M}_{n,1}=n-1$.

If we define for $k=1,\cdots,n-1$ 
$$w_k=\opo_k\, v_0^{(n\,c)}\quad \text{with}\quad v_0^{(n\,c)}=\underbrace{h_0^{(1)}\otimes h_0^{(1)} \cdots \otimes\, h_0^{(1)}}_n\ ,$$
(where all the representations are labeled by $c$ and $\gamma$) we have from (\ref{sigma0}) that $\sigma_i\;v_0^{(n\, c)}= q^{-2\,c\gamma} v_0^{(n\, c)}$. 

From (\ref{sigmaopo}), after a renormalization 
of $\sigma_i$ by the constant factor $q^{-2\,c\gamma}$, we get:
\begin{equation}\label{Burau}
\begin{split}
 \sigma_k\; w_k&=-q^{-2\gamma}\,w_k\ ,\\
 \sigma_k\; w_{k+1}&= q^{-\gamma}\,w_k+ w_{k+1}\ ,\\
 \sigma_k\; w_{k-1}&=w_{k-1}+q^{-\gamma}\, w_k\ ,\\
 \sigma_k\; w_{k+i}&=w_{k+i} \quad \text{for}\ |i|>1\ . 
 \end{split}
\end{equation}

This is the {\it reduced Burau representation} (see Section 3.3 of \cite{KT} with the rescaling on the vectors $q^{k\gamma} w_k\to b_{n-k}$ and the substitution $q^{-2\gamma}\to t$).

Next we  consider again $N=1$ but we take one representation in the tensor product different with respect to the others. 

We define $w_k^{(j)}=\opo_k (h_0^{(1)}\otimes\, \cdots\otimes\,\underbrace{h_0^{(2)}}_j\otimes\, h_0^{(1)})$ and we get the following representation for the $\sigma_i$:

\begin{equation*}\label{N=1c_1c_2}
\begin{split}
i>j+1&\quad\text{or}\quad i<j-2\\
 \sigma_i\; w_i^{(j)}&=-q^{-2\gamma_1}q^{-2c_1\gamma_1}\; w_i^{(j)}\ ,\\
 \sigma_i\; w_{i+1}^{(j)}&= q^{-\gamma_1}q^{-2c_1\gamma_1}\; w_i^{(j)}+ q^{-2c_1\gamma_1}\; w_{i+1}^{(j)}\ ,\\
 \sigma_i\; w_{i-1}^{(j)}&=q^{-2c_1\gamma_1}\; w_{i-1}^{(j)}+q^{-\gamma_1}q^{-2c_1\gamma_1}\; w_i^{(j)}\ ,\\
 \sigma_i\; w_{i+k}^{(j)}&=q^{-2c_1\gamma_1}\; w_{i+k}^{(j)}\qquad \text{for}\qquad |k|>1\ ,\\
 \\
i=j-2&\\
 \sigma_i\; w_i^{(j)}&=-q^{-2\gamma_1}q^{-2c_1\gamma_1}\; w_i^{(j)}\ ,\\
 \sigma_i\; w_{i+1}^{(j)}&= q^{-\gamma_1}q^{-2c_1\gamma_1}[\gamma_1]_q^{-1/2}[\gamma_2]_q^{1/2}\; w_i^{(j)}+ q^{-2c_1\gamma_1}\; w_{i+1}^{(j)}\ ,\\
 \sigma_i\; w_{i-1}^{(j)}&=q^{-2c_1\gamma_1}\; w_{i-1}^{(j)}+q^{-\gamma_1}q^{-2c_1\gamma_1}\; w_i^{(j)}\ ,\\
 \sigma_i\; w_{i+k}^{(j)}&=q^{-2c_1\gamma_1}\; w_{i+k}^{(j)}\qquad \text{for}\qquad |k|>1\ ,\\
 \\
i=j-1&\\
 \sigma_i\; w_i^{(j)}&=-q^{-\gamma_1-\gamma_2}q^{-c_2\gamma_1-c_1\gamma_2}\; w_i^{(j-1)}\ ,\\
 \sigma_i\; w_{i+1}^{(j)}&= q^{-\gamma_1}q^{-c_2\gamma_1-c_1\gamma_2}\; w_i^{(j-1)}+ q^{-c_2\gamma_1-c_1\gamma_2}[\gamma_1]_q^{-1/2}[\gamma_2]_q^{1/2}\; w_{i+1}^{(j-1)}\ ,\\
 \sigma_i\; w_{i-1}^{(j)}&=[\gamma_1]_q^{1/2}[\gamma_2]_q^{-1/2}\,(q^{-c_2\gamma_1-c_1\gamma_2}\;w_{i-1}^{(j-1)}+
 q^{-\gamma_2}q^{-c_2\gamma_1-c_1\gamma_2}\; w_i^{(j-1)}\,)\ ,\\
 \sigma_i\; w_{i+k}^{(j)}&=q^{-c_2\gamma_1-c_1\gamma_2}\; w_{i+k}^{(j-1)}\qquad \text{for}\qquad |k|>1\ ,\\
 \end{split}
 \end{equation*}
 \begin{equation*}\label{N=1c_1c_2-bis}
 \begin{split}
i=j\qquad&\\
 \sigma_i\; w_i^{(j)}&=-q^{-\gamma_1-\gamma_2}q^{-c_2\gamma_1-c_1\gamma_2}\; w_i^{(j+1)}\ ,\\
 \sigma_i\; w_{i+1}^{(j)}&= [\gamma_1]_q^{1/2}[\gamma_2]_q^{-1/2}\,(q^{-\gamma_2}q^{-c_2\gamma_1-c_1\gamma_2}\; w_i^{(j+1)}+ q^{-c_2\gamma_1-c_1\gamma_2}\; w_{i+1}^{(j+1)}\,)\ ,\\
 \sigma_i\; w_{i-1}^{(j)}&=q^{-c_2\gamma_1-c_1\gamma_2}[\gamma_1]_q^{-1/2}[\gamma_2]_q^{1/2}\; w_{i-1}^{(j+1)}+ q^{-\gamma_1}q^{-c_2\gamma_1-c_1\gamma_2}\; w_i^{(j+1)}\ ,\\
 \sigma_i\; w_{i+k}^{(j)}&=q^{-c_2\gamma_1-c_1\gamma_2}\; w_{i+k}^{(j+1)}\qquad \text{for}\qquad |k|>1\ ,\\
 \\
 i=j+1&\\
 \sigma_i\; w_i^{(j)}&=-q^{-2\gamma_1}q^{-2c_1\gamma_1}\; w_i^{(j)}\ ,\\
 \sigma_i\; w_{i+1}^{(j)}&=q^{-\gamma_1}q^{-2c_1\gamma_1}\; w_i^{(j)}+q^{-2c_1\gamma_1}\; w_{i+1}^{(j)}\ ,\\
 \sigma_i\; w_{i-1}^{(j)}&= q^{-2c_1\gamma_1}\; w_{i-1}^{(j)}+q^{-\gamma_1}q^{-2c_1\gamma_1}[\gamma_1]_q^{-1/2}[\gamma_2]_q^{1/2}\; w_i^{(j)}\ ,\\
 \sigma_i\; w_{i+k}^{(j)}&=q^{-2c_1\gamma_1}\; w_{i+k}^{(j)}\qquad \text{for}\qquad |k|>1\ .
\end{split}
\end{equation*}

\bigskip
In the specific case with $n=3$ we have the following matrix expressions:
\begin{equation*}
\sigma_1=
\left(
\begin{array}{cccccc}
 -d_1\, q^{-2 \gamma_1} & 0 & 0 & d_1 d_3^{-1}\, q^{-\gamma_1} & 0 & 0 \\
 0 & 0 & -d_2\, q^{-\gamma_1-\gamma_2} & 0 & 0 & d_2 d_3\, q^{-\gamma_2} \\
 0 & -d_2\, q^{-\gamma_1-\gamma_2} & 0 & 0 & d_2\, q^{-\gamma_1} & 0 \\
 0 & 0 & 0 & d_1 & 0 & 0 \\
 0 & 0 & 0 & 0 & 0 & d_2 d_3   \\
 0 & 0 & 0 & 0 & d_2 d_3^{-1} & 0   \\
\end{array}
\right)\ ,
\end{equation*}
\begin{equation*}
\sigma_2=
\left(
\begin{array}{cccccc}
 0 & d_2d_3^{-1} & 0 & 0 & 0 & 0 \\
 d_2d_3 & 0 & 0 & 0 & 0 & 0   \\
 0 & 0 & d_1 & 0 & 0 & 0 \\
 0 & d_2\, q^{-\gamma_1} & 0 & 0 & -d_2\, q^{-\gamma_1-\gamma_2} & 0 \\
 d_2 d_3\, q^{-\gamma_2} & 0 & 0 & -d_2\, q^{-\gamma_1-\gamma_2} & 0 & 0 \\
 0 & 0 & d_1 d_3^{-1}\, q^{-\gamma_1} & 0 & 0 & -d_1 q^{-2\gamma_1} \\
\end{array}
\right)\ ,
\end{equation*}

with $\quad d_1=q^{-2c_1 \gamma_1}\ $, $\quad d_2=q^{-c_2\gamma_1-c_1\gamma_2}\quad$ and $\quad d_3=[\gamma_1]_q^{1/2}[\gamma_2]_q^{-1/2}$.

\subsection{N=2}

We would like to describe the simplest example with $N=2$.

A representation of $B_n$ is obtained by the set of elements $w_{i,j}=\opo_i\,\opo_j\, (h_0^{(1)}\otimes\, \cdots\otimes\, h_0^{(1)})$ with $1\leq i\leq j \leq n-1$;
there are ${\cal M}_{n,2}=n(n-1)/2$ vectors.

The resulting formulas for the generators of $B_n$ are  (where $\{i-1,i,i+1\}\cap\{j,k\}=0$)
\begin{align}\label{LKB}
\sigma_i\; &w_{j,k} =w_{j,k}                                        &\sigma_i\; &w_{i,i} =q^{-4\gamma}w_{i,i}\nonumber\\
\sigma_i\; &w_{i,k} =-q^{-2\gamma} w_{i,k}                          &\sigma_i\; &w_{k,i}=-q^{-2\gamma} w_{k,i}\nonumber\\
\sigma_i\; &w_{i+1,k} =q^{-\gamma} w_{i,k}+w_{i+1,k}                &\sigma_i\; &w_{k,i+1}=q^{-\gamma} w_{k,i}+w_{k,i+1}\nonumber\\
\sigma_i\; &w_{i-1,k} =q^{-\gamma} w_{i,k}+w_{i-1,k}                &\sigma_i\; &w_{k,i-1}=q^{-\gamma} w_{k,i}+w_{k,i-1}\nonumber\\
\sigma_i\; &w_{i-1,i} =-q^{-2\gamma}(w_{i-1,i}+ q^{-\gamma} w_{i,i}) &\sigma_i\; &w_{i,i+1} =-q^{-2\gamma}(w_{i,i+1}+ q^{-\gamma} w_{i,i}) \nonumber\\
\sigma_i\; &w_{i-1,i-1} =q^{-2\gamma}w_{i,i}+ 2q^{-\gamma} w_{i-1,i}+w_{i-1,i-1} &&\nonumber\\
\sigma_i\; &w_{i+1,i+1} =q^{-2\gamma}w_{i,i}+ 2q^{-\gamma} w_{i,i+1}+w_{i+1,i+1}  &&\nonumber\\
\sigma_i\; &w_{i-1,i+1} =q^{-2\gamma}w_{i,i}+ q^{-\gamma} (w_{i-1,i}+w_{i,i+1})+w_{i-1,i+1}&& 
\end{align}

This is the representation given in formula (44) of \cite{JK}  using the following correspondence (with their elements $w$ rewritten as $W$ and their $q$ equal to $1$):
\begin{align}\label{corrJK}
&w_{i,i}=-2\, W_{i,i+1}\quad, \nonumber\\
&w_{i,i+1}=s^{-1}\,W_{i,i+1}-W_{i,i+2}+s\, W_{i+1,i+2}\quad,\\
&w_{i,r}=-\,W_{i,r+1}+s\,W_{i+1,r+1}+s^{-1}\, W_{i,r}-\,W_{i+1,r}\quad,\quad r\geq i+2\ .\nonumber
\end{align}

In \cite{JK} it is shown the isomorphism between the representation ${\bf W}_{n,2}$ given in their formula (44) and the Lawrence-Krammer-Bigelow representation. 

In the specific $n=3$ case we get the following matrices:
\begin{align}\label{LKB-3}
\sigma_1=\left(
\begin{array}{ccc}
 q^{-4 \gamma} & 0 & 0 \\
 -q^{-3 \gamma} & -q^{-2 \gamma} & 0 \\
 q^{-2 \gamma} & 2 q^{-\gamma} & 1 \\
\end{array}
\right)\ , &&
\sigma_2=\left(
\begin{array}{ccc}
 1 & 2 q^{-\gamma} & q^{-2 \gamma} \\
 0 & -q^{-2 \gamma} & -q^{-3 \gamma} \\
 0 & 0 & q^{-4 \gamma} \\
\end{array}
\right)\ .
\end{align}

\section{Conclusions}

In this paper we show how is possible to obtain representations of the braid groups starting from a quantum enveloping algebra whose classical limit is the Lie algebra of 
the harmonic oscillator.

The fact that the representations of this algebra are infinite dimensional obliges us to work with the lowest weights vectors that form a representation of the braid group 
due to the commutativity of the operator ${\cal P}\Re$ with the coproduct. The formulas (\ref{sigmaopo}) permit to construct a braided version of 
a classical representation of the braid group. It would be interesting to study the possible relations with the link invariants following the quantum groups techniques.


\newpage

\begin{thebibliography}{3}
 \bibitem{A}     E. Artin; {\it Theorie der Z\"opfe}, Abh. Math. Sem. Univ. Hamburg. {\bf 4}, 47--72, (1926)
 \bibitem{Burau} W. Burau; {\it \"Uber Zopfgruppen und gleichsinnig verdrillte Verkettungen}, Abh. Math. Sem. Univ. Hamburg {\bf 11}, 179--186, (1936)
 \bibitem{CP}    V. Chari, A. Pressley; {\it A guide to quantum groups}, Cambridge University Press, Cambridge (1995)
 \bibitem{LKB}   S. Bigelow; {\it The Lawrence-Krammer representation}, Topology and geometry of manifolds, Proc. Sympos. Pure Math. {\bf 71}, 51--68, (2003)
 \bibitem{3dc}   E. Celeghini, R. Giachetti, E. Sorace, M. Tarlini; {\it 3-Dimensional quantum groups from contraction of $SU(2)_q$}, J. Math. Phys. {\bf 31}, 2548--2551, (1990)
 \bibitem{h1q}   E. Celeghini, R. Giachetti, E. Sorace, M. Tarlini; {\it The quantum Heisenberg group $H(1)_q$}, J. Math. Phys. {\bf 32}, 1155--1158, (1991)
 \bibitem{FRT}   L.D. Faddeev, N. Yu. Reshetikhin, L. A. Takhtajan; {\it Quantization of Lie Groups and Lie Algebras},Leningrad Math. J., {\bf 1}, 193--225 (1989)
 \bibitem{KR}    A.N. Kirillov,  N. Yu. Reshetikhin; {\it Representations of the algebra $U_q(sl_)$, q-orthogonal polynomials and invariants of links}, 
                 in ``Infinite-dimensional Lie algebras and groups'' (Luminy-Marseille, 1988), 285--339, World Scientific Publ., Teaneck, NJ, 1989
 \bibitem{RT}    N. Yu. Reshetikhin, V. G. Turaev; {\it Ribbon graphs and their invariants derived from quantum groups}, Comm. Math. Phys. {\bf 127}, 1--26 (1990)
 \bibitem{GS}    C. Gomez, G. Sierra; {\it Quantum harmonic oscillator algebra and link invariants}, J. Math. Phys. {\bf 34}, 2119--2131, (1993) 
 \bibitem{JK}    C. Jacson, T. Kerker; {\it The Lawrence-Kramer-Bigelow representations of the braid groups via $U_q(sl_2)$}, Adv. in Math. {\bf 228}, 1689--1717, (2011) 
 \bibitem{KT}    C. Kassel, V. Turaev; ``Braid Groups'', Graduate Text in Mathematics {\bf 247}, Springer Science, New York, NY (2008) doi:10.1007/978-0-387-68548-9.
\end{thebibliography}
\end{document}